\documentclass[hidelinks,10pt]{article}

\usepackage{amsmath,amsthm,amssymb,stmaryrd,bigints}
\usepackage[pdftex]{graphicx}
\usepackage{caption,subcaption,mathabx}
\usepackage[all]{xy}
\usepackage{datetime,mathabx,color,cancel,listings}
\usepackage[hyperfootnotes=false]{hyperref}

\usepackage[margin=2.5cm]{geometry}

\newcommand{\e}{\varepsilon}

\newcommand{\R}{\mathbb{R}}

\newcommand{\beq}{\begin{equation}}
\newcommand{\ee}{\end{equation}}
\newcommand{\bac}{\begin{array}{c}}
\newcommand{\ea}{\end{array}}
\newcommand{\bal}{\begin{aligned}}
\newcommand{\eal}{\end{aligned}}

\newcommand{\real}{\operatorname{Re}}
\newcommand{\imag}{\operatorname{Im}}

\newcommand{\vertiii}[1]{{\left\vert\kern-0.25ex\left\vert\kern-0.25ex\left\vert #1 
    \right\vert\kern-0.25ex\right\vert\kern-0.25ex\right\vert}}

\begin{document}
\newtheorem{theorem}{Theorem}[section]
\newtheorem{lemma}[theorem]{Lemma}
\newtheorem{remark}[theorem]{Remark}
\newtheorem{observation}[theorem]{Observation}
\newtheorem{definition}[theorem]{Definition}
\newtheorem{example}[theorem]{Example}
\newtheorem{corollary}[theorem]{Corollary}
\newtheorem{assumption}{Assumption}
\newtheorem{property}{Property}

\title{ Semiclassical regularization of Vlasov equations\\ and wavepackets for nonlinear Schr\"odinger equations}

\author{Agissilaos Athanassoulis \footnote{agis.athanassoulis@le.ac.uk, Department of Mathematics, University of Leicester, UK}}

\maketitle

\begin{abstract}
We consider the semiclassical limit of nonlinear Schr\"odinger equations with wavepacket initial data. We recover the Wigner measure of the problem, a macroscopic phase-space  density which controls the propagation of the physical observables such as mass, energy and momentum. Wigner measures have been used to create effective models for wave propagation in random media, quantum molecular dynamics, mean field limits, and the propagation of electrons in graphene. In nonlinear settings, the Vlasov-type equations obtained for the Wigner measure are often ill-posed on the physically interesting spaces of initial data. In this paper we are able to select  the  measure-valued solution of the 1+1 dimensional Vlasov-Poisson equation which correctly captures the semiclassical limit, thus finally resolving the non-uniqueness in the seminal result of [Zhang, Zheng \& Mauser, Comm. Pure Appl. Math.
(2002) 55, doi:10.1002/cpa.3017]. The same approach is also applied to the Vlasov-Dirac-Benney equation with small wavepacket initial data, extending several known results.
\end{abstract}



\tableofcontents

\section{Introduction}

\subsection{The problem}

A well known asymptotic problem for nonlinear Schr\"odinger equations
\beq\label{eq:nls1}
i\e \partial_t \psi^\e + \frac{\e^2}2 \Delta \psi^\e - F(|\psi^\e|^2) \psi^\e =0, \qquad \psi^\e(t=0)=\psi_0^\e \in H^1(\R^n)
\ee
is to describe the evolution of {\em macroscopic observables}, such as
\beq
\begin{array}{ll}
\mbox{mass} & m(x,t) = |\psi^\e(x,t)|^2, \\
\mbox{momentum} & j(k,t) = \e^n |\widehat{\psi}^\e(\e k,t)|^2, \\
\mbox{kinetic energy} & E_{kin}(x,t) = |\nabla \psi^\e(x,t)|^2
\ea
\ee
when $\e \to 0$. Variations of this problem arises in many different physical contexts, including quantum molecular dynamics \cite{AFFPG}, mean field limits \cite{bene,GMP,GP}, wave propagation over large (geophysical) distances \cite{RKP,smit}, the formation of rogue waves \cite{Dambr} and the study of graphene \cite{fj,graphen}. We will use the term semiclassical to describe this regime \cite{BZ,Jenk,Jin}; other terms used in the literature are zero-dispersion limit \cite{Tovbis}, high frequency limit \cite{gmmp}, and geometric optics \cite{CarlesARMA,CarlesINDI}.

While direct solution of \eqref{eq:nls1} becomes more and more expensive as $\e \to 0$, it often turns out that we can recover approximations to the observables with $O(1)$ cost, i.e. with complexity independent of $\e$. This can be achieved by taking a quadratic transform of \eqref{eq:nls1}, namely the {Wigner transform} (WT)
\beq
\label{eq:WT}
W^\e(x,k,t)=
W^\e[\psi(t)](x,k) = \int\limits_{y} e^{-2\pi i k \cdot y} \psi(x+\frac{\e y}2,t) \overline{\psi}(x-\frac{\e y}2,t)dy,
\ee
leading to the Wigner equation
\beq
\label{eq:WignerEq1}
\bac
\partial_t W^\e + 2\pi k \cdot \nabla_x +i \mathcal{F}_{K\to k}^{-1} \left[
\frac{ V(x+\frac{\e K}2,t) - V(x-\frac{\e K}2,t)  }{\e} \mathcal{F}_{k'\to K} \left[ W^\e(x,k',t) \right]
\right] =0, \\ 
V(x,t) = F\left( \int\limits_\xi W^\e(x,\xi,t) \right).
\ea
\ee
This is essentially a second moment of \eqref{eq:nls1}, and it has two important properties. First of all, equation \eqref{eq:WignerEq1} has a meaningful (formal, for now) limit as $\e\to 0$, namely the Vlasov-type equation
\beq
\label{eq:VlasovEq1}
\bac
\partial_t W^0 + 2\pi k \cdot \nabla_x W^0 - \frac{1}{2\pi} \nabla_x V \cdot \nabla_k W^0  =0, \qquad
V(x,t) = F\left( \int\limits_\xi W^0(x,\xi,t) \right).
\ea
\ee
Moreover, the {\em Wigner measure}, i.e. the limit  of the Wigner transform
\beq \label{eq:WM}
W^0 = \lim\limits_{\e \to 0} W^\e 
\ee
controls macroscopic observables \cite{LP,gmmp}, e.g.
\beq
\begin{array}{ll}
\mbox{mass} & m(x,t) \approx \int\limits_k W^0(x,k,t) dk, \\
\mbox{momentum} & j(k,t) \approx \int\limits_x W^0(x,k,t) dx, \\
\mbox{kinetic energy} & E_{kin}(x,t) \approx 4\pi^2 \int\limits_k |k|^2 W^0(x,k,t) dk,
\ea
\ee
etc. A self-contained discussion of Wigner measures, including  the sense of convergence and the systematic extraction of observables, can be found in Section \ref{sec:WMs}.

This technique has been  established for a wide variety of wave problems, including Schr\"odinger \cite{AFFPG,Ath1,AP,AP2,8,AMP,bal,Bardos,Bardos2,bene,gmmp,GMP,GP,LP,pinaud1,ZZM}, Dirac \cite{fj,graphen}, and acoustic \cite{bal,miller}, elastic and Maxwell equations with smooth, random or periodic coefficients \cite{borcea,gmmp,RKP}.

A key trade-off between this approach and 
WKB-type expansions \cite{CarlesARMA,CarlesINDI,grenier,Jenk,Jin,Kamv,Tovbis} is that  we no longer try to approximate $\psi^\e$, but only the observables, through the Wigner measure.  In return, we get an elegant and widely applicable model, including in many cases the painless resolution of caustics. This can be seen as a semiclassical regularisation and continuation of the WKB system past the formation of caustics, by the introduction of a novel sense of solution \cite{JinLi}. Moreover, approximations of $\psi^\e(t)$ are often destroyed by nonlinear effects at much faster than macroscopic approximation for $W^\e[\psi^\e(t)]$; this can be seen very clearly in the discussion after Theorem \ref{thrm:Main2}.

Another important advantage of the Wigner measures approach is that it is completely non-parametric, thus being appropriate for noisy problems where the data of interest are not of WKB or other explicit parametric forms \cite{RKP,smit,Smit}. In fact, the second-moment character of the Wigner transform makes it a particularly powerful tool for stochastic problems, and it has played a key role in the recent understanding of self-averaging in wave propagation in random media \cite{bal,pinaud1}. In the same context, the Wigner transform  seems to be the appropriate generalization of the spectral  density for harmonizable (non-stationary) processes \cite{MF}.

Infinite systems of Schr\"odinger equations can be treated with Wigner measures using the same formalism; this aspect is crucial in certain fields such as statistical physics \cite{8,bene,GMP,GP}. It must be noted that infinite systems of nonlinear Schr\"odinger equations (often referred to as ``mixed states'') are attracting intense attention recently \cite{chen,LS}, following recent fundamental advances in harmonic analysis \cite{Frank}. In fact, in the context of Wigner measures, mixed states lead to simpler problems as they lead to initial data $W^0_0$ in Sobolev spaces, or even in spaces of analytic functions. This is elaborated e.g. in \cite{brezzi,LP}. In this work we will focus on pure states only, i.e. we will always start from a single nonlinear Schr\"odinger equation \eqref{eq:nls1}.

While for many classes of problems the Wigner measures approach is worked out, key questions are still open in many interesting problems, such as systems with eigenvalue crossings \cite{fj,graphen}, nonsmooth \cite{AFFPG,Ath1,AP,AP2}, and nonlinear problems. In nonlinear problems in particular, the limit Vlasov-type equation \eqref{eq:VlasovEq1} is typically not well-posed for measures. For example, in the seminal work by Zhang, Zheng \& Mauser \cite{ZZM}, it is shown that if we start with the $1$-dimensional Schr\"odinger-Poisson equation,
\beq\label{eq:nls-poiss1}
i\e \partial_t \psi^\e + \frac{\e^2}2 \Delta \psi^\e - \frac{b}2 \int\limits_{y} |x-y||\psi^\e(y,t)|^2dy \,\, \psi^\e =0, \qquad \psi^\e(t=0)=\psi_0^\e \in H^1(\R^n)
\ee
its Wigner measure $W^0 = \lim\limits_{\e\to 0} W^\e[\psi^\e]$ satisfies  (in an appropriate weak sense \cite{ZM}) the $1+1$-dimensional Vlasov-Poisson equation with initial data $W^0_0 = \lim\limits_{\e\to 0} W^\e[\psi^\e_0]$.  However, the notion of solution  used for the Vlasov-Poisson equation is so weak that uniqueness is lost. The question of determining the correct weak solution for the semiclassical limit has been the subject of numerical investigation \cite{Jnvp}, but it is still not settled. Theorem \ref{thrm:Main1} answers this question for any wavepacket initial data.

More recently, Bardos \& Besse in the breakthrough paper \cite{Bardos2} showed that, under appropriate conditions, in the case of the defocusing cubic nonlinearity
\beq\label{eq:nls-cubic1}
i\e \partial_t \psi^\e + \frac{\e^2}2 \Delta \psi^\e - b |\psi^\e|^2 \psi^\e =0, \qquad \psi^\e(t=0)=\psi_0^\e \in H^1(\R^n)
\ee
the Wigner measure indeed satisfies the resulting Vlasov-Dirac-Benney equation
\beq
\label{eq:VlasovDiracEq1}
\bac
\partial_t W^0 + 2\pi k \cdot \nabla_x + \frac{b}{2\pi} \nabla_x V \cdot \nabla_k W^0  =0, \qquad 
V(x,t) =  \int\limits_\xi W^0(x,\xi,t).
\ea
\ee
However this equation is known to be ill-posed on any Sobolev space \cite{Bardos}, and at the moment there is no sense of measure-valued solutions.

Thus, the picture that emerges for nonlinear problems can be described as follows: in many cases a Vlasov equation can be derived and justified, i.e. it can be shown that the Wigner measure does satisfy it. However this is only the first step towards approximating the evolution of the Wigner measure in time, as the Vlasov equation may be ill-posed. Indeed, as we saw, neither uniqueness nor stability can be taken for granted. In this paper  we construct an approximation to the Wigner measure for wavepackets evolving under some common nonlinearities for long times, thus extracting the correct weak solution for the semiclassical limit. It must be noted that this approach is completely non-parametric, being based only on space and Fourier localization of the initial data $\psi^\e_0$.

 The main results are stated in Section  \ref{sec:mainthrms}.  Comparisons with existing (positive and negative) results are also given.   The proofs of the main results can be found in Section \ref{sec:mainproofs}, while auxiliary results are stated and proved in Sections \ref{sec:WMs} and \ref{sec:Bckgr}.


\subsection{Notations and Definitions}\label{sec:introNot}

We will use standard multi-index notations. The Fourier transform normalization will be
\[
\widehat{f}(k) = \int\limits_{x\in \R^n} e^{-2\pi i k \cdot x} f(x) dx.
\]
Because of the particular manipulations necessary in this work, we will keep track of variable names under Fourier transforms with the notation
\begin{equation*}
\begin{aligned}
\widehat{f}(k)= \mathcal{F}_{x\to k}[f] = &\int\limits_{x\in \R^n} e^{-2\pi i k \cdot x} f(x) dx, \\
\widehat{f}(X,K) = \mathcal{F}_{x,k\to X,K}[f] = &\int\limits_{x,k\in \R^n} e^{-2\pi i [x\cdot X + k \cdot K]} f(x,k) dxdk,\\
\mathcal{F}_{k\to K}[f] = &\int\limits_{k\in \R^n} e^{-2\pi i  k \cdot K} f(x,k) dk.
\end{aligned}
\end{equation*}

The convention $\widehat{X} := \{ f \,\big|\, \widehat{f} \in X \}$ will be used for brevity.


We will use the Wiener-Sobolev spaces $A^s$. They are introduced in Definition \ref{def:new0001} in  Section \ref{sec:WMs}, along with some motivation and context. 

We will denote by $\mathcal{T}_{z}$ the translation operator
\beq
\mathcal{T}_{z} f (x) = f(x+z),
\ee
and by $\mathcal{M}_{z}$ the modulation operator
\beq
\mathcal{M}_{z} f(x) = e^{- 2 \pi i z\cdot x} f(x).
\ee

\begin{definition} \label{def:wavepackets} Let $\psi \in H^1 \cap \widehat{H}^1$ be a wavefunction with unit mass, i.e. $\|\psi\|_{L^2}=1$. We will denote
\beq
\mu_x(\psi) :=  \int\limits_x x|\psi|^2 dx, \qquad \mu_k(\psi) := { \e} \int\limits_k k|\widehat{\psi}|^2 dk,
\ee
and read $\mu_x$ as the mean position and $\mu_k$ as the mean (rescaled) momentum of the wavefunction $\psi$.
Moreover, we will
denote
\beq
\sigma^2_x(\psi) :=  \int\limits_x \left(x-\mu_x(\psi)\right)^2|\psi|^2 dx, \qquad \sigma^2_k(\psi) :=  \e^2 \int\limits_k \left(k-\frac{\mu_k(\psi)}{\e}\right)^2  |\widehat{\psi}|^2 dk,
\ee
and read $\sigma^2_x$ as the variance in position and $\sigma^2_k$ as the variance in (rescaled) momentum of the wavefunction $\psi$.
\end{definition}

The variances $\sigma^2_x(\psi)$, $\sigma^2_k(\psi)$ are the only measures of space and Fourier localization that we use to develop our non-parametric wavepacket analysis.  It can be shown that 
\beq \label{eq:wavpknot}
\sigma_x(\psi^\e)+\sigma_k(\psi^\e)=o(1)
\ee
holds for all standard classes of parametric wavepackets, such as coherent states and squeezed states, as well as less common parametric classes like chirps. In any case, this fully non-parametric notion of wavepacket through \eqref{eq:wavpknot} is quantified by Corollary \ref{cor:wavpkts}, where it is shown that 
\[
 \|W^\e[\psi^\e](x,k) - \delta\left(x-{\mu_x(\psi^\e),k-\mu_k(\psi^\e)}\right)\|_{A^{-1}} \leqslant 2\pi \Big( \sigma_x(\psi^\e)+ \sigma_k(\psi^\e) \Big).
\]
(The Banach space $A^{-1}$, specified in Definition \ref{def:new0001}, contains $\delta$-functions.) 

It must be noted that, when working on the appropriate frame of reference, the variances $\sigma^2_x(\psi)$, $\sigma^2_k(\psi)$ take a very simple form:

\begin{observation}
If a wavefunction $\psi$ is centered via a Galilean transform, i.e. if 
\beq
u=\mathcal{M}_{\frac{\mu_k(\psi)}{\e} }\mathcal{T}_{\mu_x(\psi)} \psi, 
\ee
then one readily computes
\beq
\label{eq:obssigma}
\bac
\mu_x( u) = \mu_k(u) = 0,\qquad
\sigma_x(\psi) = \sigma_x(u) = \| x u \|_{L^2_x}, \quad 
\sigma_k(\psi) = \sigma_k(u) = \frac{1}{2\pi} \| \e \nabla u   \|_{L^2}.
\ea
\ee
\end{observation}


%

The uncertainty principle \cite{up} means that we cannot make both of $\sigma_x(\psi)$, $\sigma_k(\psi)$ arbitrarily small at the same time, e.g.
\beq
\sigma_x(\psi)  \sigma_k(\psi)  \geqslant \frac{ \e \|\psi\|_{L^2(\R)}^2}{4 \pi}.
\ee
While only gaussian coherent states saturate the uncertainty principle, equation \eqref{eq:wavpknot} outlines a much broader class. Squeezed states, a class of wavepackets  generalizing coherent states, are properly introduced in Definition \ref{def:sqstates}.

\section{Statement of the main results} \label{sec:mainthrms}

\subsection{Wigner measures for wavepackets} \label{subsec:mainres}

\begin{theorem}[$1$-dimensional defocusing Schr\"odinger-Poisson equation] \label{thrm:Main1}
Let $\psi^\e(t)$ be the solution of
\beq \label{eq:1_3}
i\varepsilon\partial_t \psi^\varepsilon + \frac{\varepsilon^2}2 \Delta \psi^\varepsilon - \frac{b}2 \int\limits_y |x-y| |\psi^\varepsilon(y,t)|^{2}dy \, \psi^\varepsilon =0, \qquad \psi^\varepsilon(t=0)=\psi^\varepsilon_0 \in \mathcal{S}(\R), \qquad \|\psi^\e_0\|_{L^2}=1
\ee
for some $b>0$.
If for some $\eta>0$
\[
\sigma_x(\psi^\e_0) < \eta, \qquad \sigma_k(\psi^\e_0) < \eta,
\]
then
\[
\left\|  {W}^\varepsilon[\psi^\e(t)]-  \delta\big(x-\mu_x(\psi^\e_0) - 2\pi t \mu_k(\psi^\e_0),k-\mu_k(\psi^\e_0)\big)\right\|_{A^{-1}}  < 
2\pi (1+t) \Big[ \eta + \sqrt{\frac{b}{2\pi} \, \eta \, \,}\Big].
\]
\end{theorem}

The proof is given in Section \ref{subsec:mainproofs1}.

\smallskip

Thus the Wigner transform for any wavepacket, i.e. any initial data $\psi^\e_0$ so that $\sigma_x(\psi^\e_0)+\sigma_k(\psi^\e_0)=o(1)$, remains close to a $\delta$-function. Moreover, despite the fact that the nonlinear effects on $\psi^\e$ are of $O(1)$, the Wigner measure is not affected  by the nonlinearity. In that sense we can say that the Wigner measure satisfies the Vlasov-Poisson equation
\beq
\partial_t W^0 + 2\pi k \cdot \nabla_x W^0+ \frac{b}{4\pi} \nabla_x \int\limits_{y,\xi}  |x-y| W^0(y,\xi,t) dyd\xi \cdot \nabla_k W^0  =0, \qquad 
W^0_0=\delta(x-x_0,k-k_0),
\ee
if the nonlinear term is completely dropped, 
which is precisely what happens if we interpret it naively\footnote{Indeed, if $W^0(x,k)=\delta(x_0,k_0)$, then 
\[
\bac
\nabla_x \int\limits_{y,\xi}  |x-y| W^0(y,\xi,t) dyd\xi \cdot \nabla_k W^0= \nabla_x \int\limits_{y,\xi}  |x-y| \delta(y-x_0,\xi-k_0) dyd\xi \cdot \nabla_k \delta(x-x_0,k-k_0) =\\

=\nabla_x \int\limits_{y}  |x-y| \delta(y-x_0) dy \cdot \nabla_k \delta(x-x_0,k-k_0) =sign(x-x_0) \nabla_k \delta(x-x_0,k-k_0).
\ea
\]
Now observe that $sign(x-x_0)$ evaluated on $x_0$ is $0$; moreover $\nabla_k \delta(x-x_0,k-k_0)$ evaluated on any $(x,k)$ with $x\neq x_0$ is $0$.}.

Moreover, Theorem \ref{thrm:Main1} remains valid for a timescale  much longer than the usual $log\frac{1}\e$ Ehrenfest time-scale \cite{CarleClo}. This can be made precise for squeezed states initial data in terms of the following

\begin{corollary}[Squeezed states for the $1$-dimensional Schr\"odinger-Poisson equation]
Let
\[
\psi^\e_0 = \e^{-\frac{n\beta}2} a(\frac{x-x_0}{\e^\beta}) e^{\frac{2\pi i k_0 \cdot (x-x_0)}\e}, \qquad 0<\beta <1,
\]
be a squeezed state as in Definition \ref{def:sqstates}, and let
\beq 
i\varepsilon\partial_t \psi^\varepsilon + \frac{\varepsilon^2}2 \Delta \psi^\varepsilon - \frac{b}2 \int\limits_y |x-y| |\psi^\varepsilon(y,t)|^{2}dy \, \psi^\varepsilon =0, \qquad \psi^\varepsilon(t=0)=\psi^\varepsilon_0.
\ee
Then there exists a constant $C$ independent of  $\e$, $t$ so that
\[
\left\|  {W}^\varepsilon[\psi^\e(t)]-  \delta\big(x-x_0 - 2\pi k_0 t ,k-k_0\big)\right\|_{A^{-1}}  < 
 (1+t) C \Big( \e^{\frac{\beta}2} + \e^{\frac{1-\beta}{2}}  \Big).
\]
\end{corollary}

\bigskip

The same approach can be applied to power nonlinearities as well:

\begin{theorem}[Defocusing power nonlinearities] \label{thrm:Main2}
Let $\psi^\e(t)$ be the solution of
\beq \label{eq:1_3_2}
i\varepsilon\partial_t \psi^\varepsilon + \frac{\varepsilon^2}2 \Delta \psi^\varepsilon - b(\e) |\psi^\e|^{2\sigma} \, \psi^\varepsilon =0, \qquad \psi^\varepsilon(t=0)=\psi^\varepsilon_0 \in \mathcal{S}(\R^n), \qquad \|\psi^\e_0\|_{L^2}=1
\ee
for some $b=b(\e)>0$.
Moreover, let $C^{GN}_*$ be the sharp constant of the Gagliardo-Nirenberg inequality, see Corollary \ref{cor:GN22} for details.
If for some $\eta>0$
\[
\sigma_x(\psi^\e_0) < \eta, \qquad \sigma_k(\psi^\e_0) < \eta, \qquad \sigma_k^{\frac{n\sigma}2}(\psi_0^\e) \, \sqrt{ \frac{b(\e)}{\e^{n\sigma}}\, \frac{C^{GN}_* (2\pi)^{n\sigma-2}}{2\sigma +2} \,\,}  < \eta,
\]
then
\[
\left\|  {W}^\varepsilon[\psi^\e(t)]-  \delta\big(x-\mu_x(\psi^\e_0) - 2\pi t \mu_k(\psi^\e_0),k-\mu_k(\psi^\e_0)\big)\right\|_{A^{-1}}  < 2\pi \Big( 3 + 2t \Big) \eta .
\]
\end{theorem}

The proof is given in Section \ref{subsec:mainproofs2}.

Allowing $b(\e)=B\e^\gamma =o(1)$ and $O(1)$ initial data, is equivalent to considering small initial data and $b=B=O(1)$, through the rescaling
\[
\bac
i\varepsilon\partial_t \psi^\varepsilon + \frac{\varepsilon^2}2 \Delta \psi^\varepsilon - \e^{\gamma}B  |\psi^\e|^{2\sigma} \, \psi^\varepsilon =0 \quad \psi^\varepsilon(t=0)=\psi^\varepsilon_0 \quad \Leftrightarrow \\
\Leftrightarrow \quad 
i\varepsilon\partial_t \Psi^\varepsilon + \frac{\varepsilon^2}2 \Delta \Psi^\varepsilon - B  |\Psi^\e|^{2\sigma} \, \Psi^\varepsilon =0
\quad \Psi^\varepsilon(t=0)= \e^{\frac{\gamma}{2\sigma}} \psi^\varepsilon_0 .
\ea
\]
 Here we keep the normalization $\|\psi^\e_0\|_{L^2}=1$ so that $W^\e[\psi^\e]$ scales correctly (i.e. so that the Wigner measure exists and is not zero). 

Note that even for these  weakly-nonlinear problems, instabilities are known to appear \cite{BZ,CarlesARMA,CarlesINDI} and the semiclassical limit for wavepackets was heretofore not known.
For example, in \cite{BZ} a model of Bose-Einstein condensates is studied, namely equation \eqref{eq:1_3_2} with
\begin{equation} \label{eq:4new}
n=3, \qquad \sigma=1, \qquad  b(\varepsilon)=\varepsilon^2>0.
\end{equation}
It is shown therein that  instabilities are possible for special localized initial data. In the same setting  it has even been shown that the Wigner measure can be discontinuous in time \cite{CarlesARMA,CarlesINDI}, also pointing towards unstable behavior. All these negative results build upon initial data of the form
$\psi_0^\e = \e^{-\frac{n}2} a(\frac{x-x_0}{\e}),$
which are localized in space but not in the Fourier variable.

It is natural to ask if for some particularly convenient initial data, like coherent states,
the semiclassical limit for \eqref{eq:4new} is known.
For  coherent states, the state of the art is \cite{CarleClo}. The main result of \cite{CarleClo} can be summarized as follows: assume
\begin{equation} \label{eq:cohstateb}
|b(\varepsilon)|=O(\varepsilon^{1+\frac{n\sigma}{2}}),
\end{equation}
and  the initial wavefunction $\psi_0^\e$ is a coherent state 
\begin{equation} \label{eq:2new}
\psi^\varepsilon_0(x) = \e^{-\frac{n}4} a(\frac{x-x_0}{\sqrt{\varepsilon}})e^{\frac{ 2\pi ik_0(x-x_0)}{\varepsilon}}, \quad a \in \mathcal{S}(\R^n), \quad \|a\|_{L^2}=1, \quad x_0,k_0 \in \mathbb{R}^n.
\end{equation}
Then this parametric form is preserved, in the sense that there exists a coherent-state approximate solution of \eqref{eq:1_3_2}, 
\begin{equation} \label{eq:cohstapp}
\|\psi^\varepsilon(x,t) - \e^{-\frac{n}4} a(\frac{x-X(t)}{\sqrt{\varepsilon}},t)e^{\frac{ 2\pi  iK(t)(x-X(t))}{\varepsilon}+i\theta(t)}\|_{L^2} = o(1),
\end{equation}
where 
$
X(t),$ $K(t),$ $\theta(t),$ $a(x,t),$
satisfy simple $\e$-independent equations.
Moreover, this is valid for timescales 
\[
t=O(log \, log \frac{1}\e).
\]
A corollary of \cite{CarleClo} is that for $|b(\varepsilon)|\geqslant \varepsilon^{1+\frac{n\sigma}{2}}$ nonlinear effects on $\psi^\e(t)$ are of $O(1)$.

Equation \eqref{eq:cohstapp} provides a lot of information for the problem, but at the cost of a rather weak nonlinearity,  i.e. assumption \eqref{eq:cohstateb}, excluding many physically relevant problems. In particular, the nonlinearity \eqref{eq:4new} is too strong for the result of \cite{CarleClo}. 
Moreover, in most realistic settings the values of $\e$ range between $10^{-2}$ and $10^{-6}$, so this would lead to short timescales as well since, for the natural logarithm, $log \, log 10^6 \approx 2.6$. 

In other words, it was non known heretofore whether we can have any control of the observables 
in the problem described by the scaling \eqref{eq:4new} for wavepacket initial data; not even for coherent state initial data. To answer this question one observes that Theorem \ref{thrm:Main2} implies the following

\begin{corollary}[Squeezed states for defocusing power nonlinearities] \label{cor:main32}
Let
\[
\psi^\e_0 = \e^{-\frac{n\beta}2} a(\frac{x-x_0}{\e^\beta}) e^{\frac{2\pi i k_0 \cdot (x-x_0)}\e}, \qquad 0<\beta <1,
\]
be a squeezed state as in Definition \ref{def:sqstates}, and let
\beq 
i\varepsilon\partial_t \psi^\varepsilon + \frac{\varepsilon^2}2 \Delta \psi^\varepsilon - \e^\gamma  |\psi^\e|^{2\sigma} \, \psi^\varepsilon =0,
\qquad \psi^\varepsilon(t=0)=\psi^\varepsilon_0.
\ee
Then there exists a constant $C$ independent of   $\e$, $t$ so that
\[
\left\|  {W}^\varepsilon[\psi^\e(t)]-  \delta\big(x-x_0 - 2\pi k_0 t ,k-k_0\big)\right\|_{A^{-1}}  < 
 (1+t) C \Big( \e^{{\beta}} + \e^{{1-\beta}} + \e^{\frac{\gamma-n\sigma \beta}{2}} \Big).
\]
\end{corollary}

Setting $\gamma=2$, $n=3$ in Corollary \ref{cor:main32} above means we recover the setting of \eqref{eq:4new}. Then, if $\psi^\e_0$ is a squeezed state with $\beta < \frac{2}3$ it follows that $W^\e[\psi^\e(t)]$ evolves linearly as long as $ t \cdot  ( \e^{{\beta}}  + \e^{1-\frac{3 \beta}{2}} )=o(1)$.

\bigskip

We can apply this approach to focusing power nonlinearities as well:

\begin{theorem}[Focusing power nonlinearities] \label{thrm:Main3}
Let $\psi^\e(t)$ be the solution of
\beq \label{eq:1_3_23}
i\varepsilon\partial_t \psi^\varepsilon + \frac{\varepsilon^2}2 \Delta \psi^\varepsilon - b(\e) |\psi^\e|^{2\sigma} \, \psi^\varepsilon =0, 
\qquad \psi^\varepsilon(t=0)=\psi^\varepsilon_0 \in \mathcal{S}(\R^n), \qquad \|\psi^\e_0\|_{L^2}=1
\ee
for some $b=b(\e)<0$, and for $n\sigma=1$.
If for some $\eta>0$
\[
\sigma_x(\psi^\e_0) < \eta, 
\qquad \sigma_k(\psi^\e_0) < \eta, 
\qquad
\frac{|b(\e)|}{\e} 
\frac{ C^{GN}_*}{ 2\pi (4+\frac{4}n)} +  \sqrt{   \frac{|b(\e)|}{\e} 
\frac{ C^{GN}_*}{ 2\pi (4+\frac{4}n)}}
 \sqrt{  \frac{|b(\e)|}{\e} 
\frac{ C^{GN}_*}{ 2\pi (4+\frac{4}n)}  + \frac{\sigma_k(\psi^\e_0)}{2\pi}
 }
 < \eta,
\]
then
\[
\left\|  {W}^\varepsilon[\psi^\e(t)]-  \delta\big(x-\mu_x(\psi^\e_0) - 2\pi t \mu_k(\psi^\e_0),k-\mu_k(\psi^\e_0)\big)\right\|_{A^{-1}}  < 2\pi \Big( 3 + 2t \Big) \eta .
\]
\end{theorem}

The proof is given in Section \ref{subsec:mainproofs3}.
The restriction $n\sigma=1$ has to do with working out explicitly the upper bound in the technical Lemma \ref{lm:bootsrmine}.

The aforementioned result of \cite{CarleClo} applies in the same way to focusing and defocusing problems. Theorem \ref{thrm:Main3}  allows for stronger focusing nonlinearities, longer timescales, and of course more general initial data. This can be seen clearly in the following

\begin{corollary}[Squeezed states for focusing nonlinearities]
Let
\[
\psi^\e_0 = \e^{-\frac{\beta}2} a(\frac{x-x_0}{\e^\beta}) e^{\frac{2\pi i k_0 \cdot (x-x_0)}\e}, \qquad 0<\beta <1,
\]
be a squeezed state as in Definition \ref{def:sqstates}, $n\sigma=1$, and 
\beq 
i\varepsilon\partial_t \psi^\varepsilon + \frac{\varepsilon^2}2 \Delta \psi^\varepsilon + \e^\gamma |\psi^\e|^{2\sigma} \, \psi^\varepsilon =0, 
 \qquad \psi^\varepsilon(t=0)=\psi^\varepsilon_0.
\ee
Then there exists a constant $C$ independent of $\e$, $t$ so that
\[
\left\|  {W}^\varepsilon[\psi^\e(t)]-  \delta\big(x-x_0 - 2\pi k_0 t,k-k_0\big)\right\|_{A^{-1}}  < 
 (1+t) C \Big( \e^{{\beta}} + e^{{1-\beta}}+ \e^{\gamma-1} + \e^{ \frac{\gamma-\beta}{2}}  \Big)
\]
\end{corollary}
Thus, for any $\gamma>1$ control of the Wigner measure is obtained, as opposed to $\gamma>\frac{3}2$ in \cite{CarleClo}.

\subsection{Idea of the proofs}

The idea behind the proofs for all of the main results follows the same general steps, bringing together several different ideas, and adjusting the details as needed for each nonlinearity: 

\noindent {\bf Step 1: Go to the appropriate frame of reference.} The nonlinearities we work with are Galilean invariant. In that context, we use a frame of reference that centers the initial data
\beq
u_0^\e(x)=\mathcal{M}_{\frac{\mu_k(\psi_0^\e)}{\e} }\mathcal{T}_{\mu_x(\psi_0^\e)} \psi_0^\e=\psi_0^\e(x+x_0)e^{- 2\pi i\frac{\mu_k(\psi_0^\e)\cdot x}{\varepsilon}}, 
\ee
and work on problem \eqref{eq:nls1} through 
\beq\label{eq:nls1GalTra}
i\e \partial_t u^\e + \frac{\e^2}2 \Delta u^\e - F(|u^\e|^2) u^\e =0, \qquad u^\e(t=0)=u_0^\e .
\ee

The Galilean invariance of \eqref{eq:nls1} (recalled in Lemmata  \ref{lm:GlInvII}, \ref{lm:GlInv}) means that $\psi^\e(x,t)$ is related to $u^\e(x,t)$ through
\[
\psi^\e(x,t) = u^\e(x-vt-x_0,t) e^{i\left( \frac{v \cdot (x-x_0)}{\e}  - \frac{v\cdot v}{2\e}\right)}, \qquad v = 2\pi \mu_k(\psi_0^\e), \qquad x_0 = \mu_x(\psi_0^\e).
\]

\noindent {\bf Step 2: Show that if $\sigma_x(\psi_0^\e),$ $\sigma_k(\psi_0^\e)$ are small, then $\sigma_x(u^\e(t)),$ $\sigma_k(u^\e(t))$ are also small.}
By state of the art methods for nonlinear Schr\"odinger equations, one can obtain bounds for $\|\e\nabla u^\e(t)\|_{L^2}$ in terms of $\|\e\nabla u^\e_0\|_{L^2}$. 
Then we proceed to bound $\|x u^\e(t)\|_{L^2}$  by appropriate functions of $\|x u^\e_0\|_{L^2}$, $\|\e \nabla u^\e_0\|_{L^2}.$ 
From this we conclude that  $\sigma_x(u^\e(t)),$ $\sigma_k(u^\e(t))$ are bounded by appropriate functions of $\sigma_x(u^\e_0)=\sigma_x(\psi^\e_0),$ $\sigma_k(u^\e_0)=\sigma_k(\psi^\e_0)$.

Working out the details in each case determines the constants and, crucially, the timescales for which this bound is useful.

\noindent {\bf Step 3: Conclude that $W^\e[u^\e(t)] \approx \delta(x-0,k-0),$ quantify the rate and timescale of convergence, and go back to the initial frame of reference to obtain the result for $W^\e[\psi^\e(t)]$.}
The previous step is exploited through Corollary \ref{cor:wavpkts} to complete the proof.

\smallskip

Every effort has been made to state and prove
regularity results, bootstrap arguments etc in a self-contained way in Sections \ref{sec:WMs} and \ref{sec:Bckgr}. That way Section \ref{sec:mainproofs} is devoted to presenting coherently  how the different pieces fit together, without being sidetracked by various technical details. The engine behind the proofs is Lemma \ref{lm:loc11} and its Corollary \ref{cor:wavpkts}, which translate $H^1$ and $\widehat{H}^1$ estimates to convergence results for the Wigner measure. It is through Lemma \ref{lm:loc11} that the new functional framework, introduced in detail in Section \ref{sec:WMs} below, makes the results of this paper possible.

\section{Wigner measures and the new functional framework} \label{sec:WMs}

The Wigner transform (WT)  can be seen as a sesquilinear transform
\[
W^\varepsilon:L^2(\R^n) \times L^2(\R^n) \to L^2(\R^{2n}): f,g \mapsto W^\varepsilon[f,g],
\] 
defined as 
\beq
\label{eq:wt}
W^\varepsilon[f,g](x,k)=\int\limits_{y \in \R^n} e^{-2\pi i k y}f(x+\frac{\varepsilon y}2)\bar g(x-\frac{\varepsilon y}2)dy.
\ee
One easily checks the following elementary properties \cite{AP3,AMP}:
\beq
\label{eq:wtspaces}
\begin{array}{r c l}
f,g \in L^2(\R^{n}) & \Rightarrow & W^\e[f,g] \in L^2(\R^{2n}) \cap L^\infty(\R^{2n}), \\
f,g \in H^1(\R^{n}) \cap \widehat{H}^1(\R^{n}) & \Rightarrow & W^\e[f,g] \in H^1(\R^{2n}) \cap \widehat{H}^1(\R^{2n}), \\
f,g \in \mathcal{S}(\R^n) & \Rightarrow & W^\e[f,g] \in \mathcal{S}(\R^{2n}).
\ea
\ee

Often the quadratic version is used, in which case we denote 
\[
W^\varepsilon[f] := W^\varepsilon[f,f]. 
\]
The WT  $W^\varepsilon[f]$ describes the quadratic observables of $f$ through
\[
\int\limits_{x,k \in \R^n} W^\varepsilon[f](x,k)\,\phi(x,k)\, dxdk = \int\limits_{x\in \R^n} \overline{f(x)} \, \phi(x,\varepsilon \nabla_x)f(x)\, dx
\]
where $\phi(x,\varepsilon \nabla_x)$ is the Weyl pseudodifferential operator with symbol $\phi(x,k)$ \cite{gmmp,LP}. Thus weak approximations of $W^\varepsilon[f]$ can provide information for the quadratic observables of $f$ -- but not for its point values.

The most fruitful application of the $\varepsilon$-dependent WT is to an $\varepsilon$-dependent family of functions, $\{\psi^\varepsilon\}_\e$.
Under appropriate conditions, it is known that $W^\varepsilon[\psi^\varepsilon]$ converges in weak-$*$ sense to a probability measure $W^0$ on $\R^{2n}$ as $\varepsilon\to 0$ \cite{LP}; $W^0$ is then called the Wigner measure (WM) of the family of functions $\{\psi^\varepsilon\}_\e$. Intuitively, the WM keeps track of the limits of the observables of $\psi^\varepsilon$ as $\varepsilon \to 0$ through
\[
\lim\limits_{\varepsilon \to 0} \int\limits_{x\in \R^n} \overline{\psi^\varepsilon(x)} \, \phi(x,\varepsilon \nabla_x)\psi^\varepsilon(x)\, dx = \int\limits_{x,k \in \R^n} W^0(x,k)\,\phi(x,k)\, dxdk
\] 
while the family $\{\psi^\varepsilon\}_\e$ itself has no meaningful limit (typically $\lim\limits_{\varepsilon\to 0}\psi^\varepsilon=0$ in the sense of distributions).


The framework developed in \cite{LP} for the weak-$*$ convergence of the WT towards the WM is based on the algebra of test functions $\mathcal{A}$, generated by the norm $\|\phi\|_{\mathcal{A}}:=  \| \mathcal{F}_{k\to K}  [\phi](x,K) \|_{L^1_KL^\infty_x}$. A back-of-the-envelope calculation explains the selection of this norm in the following sense: Let $\|\psi^\varepsilon\|_{L^2}=1$, then
\begin{eqnarray}
\nonumber
\int\limits_{x,k \in \R^n} W^\varepsilon[\psi^\varepsilon](x,k) \phi(x,k) dxdk = \int\limits_{x,k,y \in \R^n} e^{-2\pi i k y} \psi^\varepsilon(x+\frac{\varepsilon y}2) \overline{ \psi^\varepsilon}(x-\frac{\varepsilon y}2) \phi(x,k) dxdk = \\
\nonumber
=\int\limits_{x,y \in \R^n} \psi^\varepsilon(x+\frac{\varepsilon y}2) \overline{ \psi^\varepsilon}(x-\frac{\varepsilon y}2) \int\limits_{k\in \R^n} e^{-2\pi i k y}  \phi(x,k) dk \,\, dxdy \quad \Rightarrow \qquad \qquad \\
\Rightarrow \quad \left| \langle W^\varepsilon[\psi^\varepsilon], \phi \rangle  \right| \leqslant \| \psi^\varepsilon(x+\frac{\varepsilon y}2) \overline{ \psi^\varepsilon}(x-\frac{\varepsilon y}2)\|_{L^\infty_y L^1_x} \| \mathcal{F}_{k\to y}[\phi] \|_{L^1_y L^\infty_x}, \qquad \qquad \label{eq:bote}
\end{eqnarray}
%
%
where of course
\[
\| \psi^\varepsilon(x+\frac{\varepsilon y}2) \overline{ \psi^\varepsilon}(x-\frac{\varepsilon y}2)\|_{L^\infty_y L^1_x} = \sup_y \int\limits_{x\in \R^n} \left| \psi^\varepsilon(x+y) \overline{ \psi^\varepsilon}(x-y) \right| dx = 1.
\]
Thus the set $\{ W^\e[\psi^\e]\}_\e$ is uniformly bounded in the dual of $\mathcal{A}$, $\mathcal{A}'$, and hence weak-$*$ compact by virtue of the Banach-Alaoglou Theorem. By extracting a subsequence in $\e$ if necessary, the WM $W^0$ is now well defined.  It is known that $W^0$ is in fact a non-negative finite measure \cite{LP}, hence the term Wigner measure is justified. 

Finding ways to metrise the weak-$*$ limit 
$$\langle W^0, \phi \rangle = \lim\limits_{\e \to 0} \langle W^\e, \phi \rangle \qquad \forall \phi \in \mathcal{A}$$
is important in itself, as it could yield better control on uniqueness questions, and of course help quantify the rate of convergence.
One might think that since $W^0$ is a probability measure,  $W^\varepsilon$ would naturally be seen converge to $W^0$ in some Banach space of measures. However, for $\psi^\varepsilon \in L^2(\R^n)$, $W^\varepsilon=W^\varepsilon[\psi^\varepsilon]\in L^2(\R^{2n})\cap L^\infty(\R^{2n})$ may not even be in $L^1(\R^{2n})$ \cite{Simon}. In that case, $\int\limits W^\varepsilon dxdk=\|\psi^\varepsilon\|_{L^2}^2$ in Cauchy-principal-value sense, but $W^\varepsilon$ does not define a finite measure at all. By using a Fourier based norm, as we do below, we go around this integrability question, and let the Fourier transform absorb any improper integrals.

\begin{definition}[The Wiener-Sobolev spaces $A^s$] \label{def:new0001}
For $s\geqslant0$, we will denote with $A^s(\R^n)$ the Banach space of functions generated by the norm
\[
\|\phi\|_{A^s(\R^n)} :=  \int\limits_{y \in \R^n} (1+|y|^s) |\widehat{\phi}(y)| dy.
\]
In phase-space this becomes
\[
\|\phi\|_{A^s(\R^{2n})} :=  \int\limits_{X,K \in \R^n} \left(1+\sqrt{|X|^2+|K|^2}^s \right) |\widehat{\phi}(X,K)| dy.
\]
When $s>0$, we will denote the dual of $A^s$ by $A^{-s}$, i.e. 
\[
\|\phi\|_{A^{-s}} = \sup\limits_{\|\psi\|_{A^s}=1} \left|  \langle \phi , \psi \rangle \right|.
\]
\end{definition}

\begin{remark} When $s=0$ we recover the standard Wiener algebra, $\|\phi\|_{A^0}=\|\widehat{\phi}\|_{L^1}$. Its dual space will be denoted as
\[
(A^0)'=\widehat{L}^\infty = \{ f: \|\widehat{f}\|_{L^\infty}< \infty \}.
\]
\end{remark}

\begin{lemma}[Consistency of $\mathcal{A}$, $A^0$ and $A^1$] \label{lm:AandclaA} For every $\phi$ in the Schwarz class of test functions $\mathcal{S}(\R^n)$
\[
\|\phi\|_{\mathcal{A}} \leqslant \| \phi \|_{A^0} \leqslant \| \phi \|_{A^1}.
\]
\end{lemma}

\noindent {\bf Proof: } Simply observe that, for any $\phi\in \mathcal{S}(\R^n)$,
\[
\bac
\|\phi\|_{\mathcal{A}}=\| \mathcal{F}_{k\to K} [\phi]\|_{L^1_K L^\infty_x}=
\int\limits_{K} \sup\limits_{x\in \R^n} |\mathcal{F}_{k\to K}[\phi](x,K)| dK 
\leqslant \int\limits_{K
}  \int\limits_{X
} \left|\mathcal{F}_{x,k\to X,K}[\phi](X,K) \right| dX  dK =\|\widehat{\phi}\|_{L^1}=\|\phi\|_{A^0}.
\ea
\]

\qed

This leads to the following

\begin{lemma} \label{thrm:WMA_1}
For any $\|\psi^\varepsilon\|_{L^2}=1$,
\[
\|W^\varepsilon[\psi^\varepsilon]\|_{A^{-1}}= \|W^\varepsilon[\psi^\varepsilon]\|_{\widehat{L}^\infty}=1.
\]
\end{lemma}

\noindent {\bf Proof:} 
First of all, recall that $\widehat{L}^\infty=(A^0)'$.
Now simply repeat the computation of equation \eqref{eq:bote}; this shows $\|W^\varepsilon[\psi^\varepsilon]\|_{\mathcal{F}L^\infty}\leqslant 1$; equality follows by selecting ${\phi}_R=e^{-\pi R(x^2+k^2)}$, and taking $\sup\limits_{R\to 0} \left| \langle W^\varepsilon[\psi^\varepsilon],\phi_R \rangle \right|$ (observe that $\|\phi_R\|_{A^0}=1$). 

The estimate $\|W^\varepsilon[\psi^\varepsilon]\|_{A^{-1}}\leqslant 1$ follows in the same way. To show that $\|W^\varepsilon[\psi^\varepsilon]\|_{A^{-1}}=1$ it suffices to take $\phi_R$ as before, and compute $\|\phi_R\|_{A^1}=1+CR^{\frac{3n}2}$.
\qed

In other words, the norms $A^{-1}$, $\widehat{L}^\infty$ are correctly scaled to capture the Wigner measure as $\e \to 0$. We will be working mainly in $A^{-1}$, that is the admissible observables will be those operators with Weyl symbols $\phi \in A^1$. Technically, this is a slightly smaller class of observables than the class $\mathcal{A}$ introduced in \cite{LP}.


%
%
%

\section{Background results} \label{sec:Bckgr}

\subsection{Background on Schr\"odinger equations}

\subsubsection{Well-posedness and conservation of energy}

The $1$-dimensional Schr\"odinger-Poisson problem has certain special features. One is that $1$-dimensional Poisson kernel, $|x|$, grows at infinity. This means that the standard methods for $V(x,t) =  \int\limits_y K(x-y) |\psi^\varepsilon(y,t)|^{2}dy$ with kernels $K \in L^\infty + L^p$ \cite{CazB} cannot be used off-the-shelf. Because of that feature, the nonlinear potential 
\beq \label{eq:1_3c}
V(x,t) = \frac{b}2 \int\limits_y |x-y| |\psi^\varepsilon(y,t)|^{2}dy
\ee
 has nontrivial behavior at infinity, 
\[
\lim\limits_{x \to \pm \infty} \frac{d}{dx} V(x,t) = \mp \frac{b}2 \|\psi^\e(x,t)\|_{L^2}^2.
\]
 We will use  the approach of \cite{ZZM}, and modify it to 
also control the moments of the solution:

\begin{theorem}[Well-posedness for the $1$-dimensional Schr\"odinger-Poisson equation] \label{thrm:wellp3}
Consider the Cauchy problem
\beq \label{eq:1_3original}
i\varepsilon\partial_t \psi^\varepsilon + \frac{\varepsilon^2}2 \Delta \psi^\varepsilon - \frac{b}2 \int\limits_y |x-y| |\psi^\varepsilon(y,t)|^{2}dy \, \psi^\varepsilon =0, \qquad \psi^\varepsilon(t=0)=\psi^\varepsilon_0 \in H^1(\R) \cap \widehat{H}^1(\R).
\ee
This problem has a unique, global-in-time solution in $H^1(\R) \cap \widehat{H}^1(\R)$ which conserves mass
\beq
\|\psi^\e(t)\|_{L^2} = \|\psi^\e_0\|_{L^2}
\ee
and energy
\beq \label{eq:spconoe}
\frac{\varepsilon^2}2 \|\nabla \psi^\varepsilon(t)\|_{L^2}^2 + \frac{b}{4} \int\limits_{x,y} |x-y| |\psi^\e(x,t)|^2 |\psi^\e(y,t)|^2 dxdy = 
\frac{\varepsilon^2}2 \|\nabla \psi^\varepsilon_0\|_{L^2}^2 + \frac{b}{4} \int\limits_{x,y} |x-y| |\psi^\e_0(x)|^2 |\psi^\e_0(y)|^2 dxdy.
\ee

Moreover,
\beq \label{eq:ststoop}
\|\e \frac{d}{dx} \psi^\e(t) \|_{L^2} \leqslant \|\e \frac{d}{dx} \psi^\e_0 \|_{L^2} + |b| \|\psi^\e_0\|_{L^2}^3 |t|  
\ee
and
\beq \label{eq:momscrpois}
\| x \psi^\e(t)\|_{L^2}\leqslant  \| x \psi^\e_0\|_{L^2} + \int\limits_{\tau=0}^t \|\e \frac{d}{dx} \psi^\e(\tau) \|_{L^2} d\tau
.
\ee
\end{theorem}

\noindent 
{\bf Proof:} By the symmetry of the problem, we readily have 
\[
\frac{d}{dt} \|\psi^\e(t)\|_{L^2}=0.
\]

Denote for brevity $V(x,t)$ the nonlinear potential as in equation \eqref{eq:1_3c}. $V(x,t)$ is the solution of 
\beq \label{eq:1_3b}
\bac
\Delta V(x,t) = b |\psi^\e(x,t)|;
\ea
\ee
either equation \eqref{eq:1_3c} or \eqref{eq:1_3b} yields
\beq
\frac{d}{dx} V(x,t) = \frac{b}2  \left( \int\limits_{y=x}^\infty |\psi^\e(y,t)|^2 dy - \int\limits_{y=-\infty}^x |\psi^\e(y,t)|^2dy \right) ,
\ee
and therefore, using the conservation of mass,
\beq
|\frac{d}{dx} V(x,t)| \leqslant \frac{|b|}2 \|\psi^\e(t)\|_{L^2}^2=\frac{|b|}2 \|\psi^\e_0\|_{L^2}^2.
\ee

Now, following  the steps of the proof
of Lemma 2.1 of \cite{ZZM}, we check that
\beq \label{eq:stepstomkjk}
\bac
\frac{1}2 \frac{d}{dt} \| \e \frac{d}{dx} \psi^\e(t) \|^2_{L^2} = -\e \imag \left[ \langle \psi^\e \frac{d}{dx} V, \frac{d}{dx} \psi^\e \rangle \right] \leqslant \frac{|b|}2 \|\psi^\e_0\|_{L^2}^3 \|\e\frac{d}{dx} \psi^\e(t)\|_{L^2} \quad \Rightarrow\\
\Rightarrow \quad \frac{d}{dt} \| \e \frac{d}{dx} \psi^\e(t) \|_{L^2} \leqslant 
|b| \|\psi^\e_0\|_{L^2}^3.
\ea
\ee
Thus equation \eqref{eq:ststoop},
which is essentially equation (2.8) of \cite{ZZM}, follows. Observe that the sign of $b$ in fact makes no difference (in \cite{ZZM} the proof is carried out for $b=1$ only).

Similarly,
\beq 
\bac
 \frac{1}2 \frac{d}{dt} \| x \psi^\e(t)\|_{L^2}^2
= \real \left[ \frac{i\e}2  \langle x \Delta \psi^\e(t), x \psi^\e(t)
 \right] 
= \real \left[ i\e  \langle  x  \psi^\e(t), \frac{d}{dx} \psi^\e(t)
 \right]
\leqslant  \| \e \nabla \psi^\e(t)\|_{L^2} \| x \psi^\e\|_{L^2} \quad \Rightarrow\\

\Rightarrow \quad 
 \frac{d}{dt} \| x \psi^\e(t)\|_{L^2}
\leqslant  \| \e \nabla \psi^\e(t)\|_{L^2}. 
\ea
\ee
Equation \eqref{eq:momscrpois} follows.

Now there is enough regularity to justify uniqueness and the conservation of energy by standard arguments \cite{CazB}.

Observe that  Lemma 2.1 of \cite{ZZM} implies $\psi^\e(t) \in H^m$ for any $m\in \mathbb{N}$ if there is sufficient regularity in the initial data. By standard arguments \cite{CazB} it follows that if there is sufficient regularity in the initial data $\|\psi^\e(t)\|_{H^1}$, $\|\psi^\e(t)\|_{\widehat{H}^1}$ are continuous functions of time.
\qed

\smallskip

Well-posedness for the nonlinear Schr\"odinger equation with power nonlinearities on $H^1$ is exhaustively well studied \cite{CazB}. Here we briefly recall the relevant results in the semiclassical scaling, and outline how control of moments ($\widehat{H}^1$ norm) follows.

\begin{theorem}[Well-posedness for energy sub-critical defocusing power nonlinearities] \label{thrm:wellp1}
Consider the Cauchy problem
\beq \label{eq:1_1}
i\varepsilon\partial_t \psi^\varepsilon + \frac{\varepsilon^2}2 \Delta \psi^\varepsilon - b |\psi^\varepsilon|^{2\sigma}\psi^\varepsilon =0, \qquad \psi^\varepsilon(t=0)=\psi^\varepsilon_0 \in H^1(\R^n)
\ee
with
\beq
b>0, \qquad 0 < \sigma < \frac{2}{(n-2)_+}.
\ee
This problem has a unique, global-in-time solution in $H^1$ which conserves mass
\beq
\|\psi^\e(t)\|_{L^2} = \|\psi^\e_0\|_{L^2}
\ee
and energy
\beq \label{eq:defocenpn}
\frac{\varepsilon^2}2 \|\nabla \psi^\varepsilon(t)\|_{L^2}^2 + \frac{b}{\sigma +1} \|\psi^\varepsilon(t)\|^{2\sigma+2}_{L^{2\sigma+2}}=
\frac{\varepsilon^2}2 \|\nabla \psi^\varepsilon_0\|_{L^2}^2 + \frac{b}{\sigma +1} \|\psi^\varepsilon_0\|^{2\sigma+2}_{L^{2\sigma+2}}
\ee
\end{theorem}

\noindent 
{\bf Proof:} The proof follows by a straightforward adaptation to the semiclassical scaling of the proof of Theorem 4.8.1 of \cite{CazB}. The result stays true if $b=b(\e)>0$.

Moreover by standard arguments \cite{CazB} it follows that if there is sufficient regularity in the initial data $\|\psi^\e(t)\|_{H^1}$, $\|\psi^\e(t)\|_{\widehat{H}^1}$ are continuous functions of time.
\qed

\begin{theorem}[Well-posedness for mass sub-critical focusing power nonlinearities] \label{thrm:wellp2}
Consider the Cauchy problem
\beq \label{eq:1_2}
i\varepsilon\partial_t \psi^\varepsilon + \frac{\varepsilon^2}2 \Delta \psi^\varepsilon - b |\psi^\varepsilon|^{2\sigma}\psi^\varepsilon =0, \qquad \psi^\varepsilon(t=0)=\psi^\varepsilon_0 \in H^1(\R^n)
\ee
with
\beq
b<0, \qquad 0 < \sigma < \frac{2}{n}.
\ee
This problem has a unique, global-in-time solution in $H^1$ which conserves mass
\beq
\|\psi^\e(t)\|_{L^2} = \|\psi^\e_0\|_{L^2}
\ee
and energy
\beq
\frac{\varepsilon^2}2 \|\nabla \psi^\varepsilon(t)\|_{L^2}^2 + \frac{b}{\sigma +1} \|\psi^\varepsilon(t)\|^{2\sigma+2}_{L^{2\sigma+2}}=
\frac{\varepsilon^2}2 \|\nabla \psi^\varepsilon_0\|_{L^2}^2 + \frac{b}{\sigma +1} \|\psi^\varepsilon_0\|^{2\sigma+2}_{L^{2\sigma+2}}
\ee
\end{theorem}

\noindent 
{\bf Proof:} The proof follows by a straightforward adaptation to the semiclassical scaling of the proof of Theorem 4.8.1 of \cite{CazB}. The result stays true if $b=b(\e)<0$.

Moreover by standard arguments \cite{CazB} it follows that if there is sufficient regularity in the initial data $\|\psi^\e(t)\|_{H^1}$, $\|\psi^\e(t)\|_{\widehat{H}^1}$ are continuous functions of time.
\qed

\begin{theorem}[Moments under power nonlinearities] \label{thrm:momscrpois22}
Let $\psi^\e$ be the solution of
\beq \label{eq:1_333}
i\varepsilon\partial_t \psi^\varepsilon + \frac{\varepsilon^2}2 \Delta \psi^\varepsilon - b |\psi^\varepsilon|^{2\sigma}\psi^\varepsilon =0, \qquad \psi^\varepsilon(t=0)=\psi^\varepsilon_0 \in H^1(\R^n).
\ee

Then
\beq \label{eq:momscrpois22}
\| x \psi^\e(t)\|_{L^2} \leqslant \| x \psi^\e_0\|_{L^2} + \int\limits_{\tau=0}^t \|\e \frac{d}{dx} \psi^\e(\tau) \|_{L^2} d\tau.
\ee
\end{theorem}

\noindent {\bf Proof:} This follows in exactly the same way as in Theorem \ref{thrm:wellp3}.
More specifically, one directly computes
\beq 
\bac
 \frac{1}2 \frac{d}{dt} \| x \psi^\e(t)\|_{L^2}^2
= \real \left[ \frac{i\e}2  \langle x \Delta \psi^\e(t), x \psi^\e(t)
 \right] 
= \real \left[ i\e  \langle  x  \psi^\e(t), \frac{d}{dx} \psi^\e(t)
 \right]
\leqslant  \| \e \nabla \psi^\e(t)\|_{L^2} \| x \psi^\e\|_{L^2} \quad \Rightarrow\\

\Rightarrow \quad 
 \frac{d}{dt} \| x \psi^\e(t)\|_{L^2}
\leqslant  \| \e \nabla \psi^\e(t)\|_{L^2}.
\ea
\ee
The result follows. 
\qed

\subsubsection{Galilean invariance}

\begin{lemma}[Galilean invariance] \label{lm:GlInvII} Let $\psi$ satisfy
\beq \label{eq:GLIoreqb}
i\varepsilon\partial_t \psi + \frac{\varepsilon^2}2 \Delta \psi - b \int\limits_y K(x-y) |\psi(y,t)|^{2} dy \, \psi=0, \qquad \psi(t=0)=\psi_0 \in L^2(\R^n).
\ee
For any $x_0,v\in \R^n$,
and denote
\beq \label{eq:GalTransfb}
u(x,t) = \psi(x+vt+x_0,t) e^{-i\left(\frac{v\cdot x}\varepsilon+\frac{v\cdot v}{2\varepsilon}t \right)}.
\ee
Then $u$ satisfies
\beq \label{eq:aljhb}
i\varepsilon\partial_t u + \frac{\varepsilon^2}2 \Delta u - b|u|^{2\sigma} u=0, \qquad u(t=0)=u_0=\psi_0(x+x_0)e^{-i\frac{v\cdot x}{\varepsilon}} \in L^2(\R^n).
\ee
Moreover,
\beq \label{eq:WTGalTransfb}
W^\e[u(t)](x,k) = W^\e[\psi(t)]\left(x+vt+x_0,k+\frac{v}{2\pi}\right).
\ee
\end{lemma}

\noindent 
{\bf Proof:} See \cite{GinibreNotes} for the transformation of equation \eqref{eq:GLIoreq}, i.e. for equations \eqref{eq:GalTransf}, \eqref{eq:aljh}.

Equation \eqref{eq:WTGalTransf} follows by the elementary computation
\[
\bac
W^\e[u(t)] = W^\e[\psi(x+vt+x_0,t) e^{-i\left(\frac{v\cdot x}\varepsilon+\frac{v\cdot v}{2\varepsilon}t \right)}] = \\

=\int\limits_y e^{-2\pi i k \cdot y} \psi(x+\frac{\e y}2 + v t + x_0,t) 
e^{-i\left(\frac{v\cdot (x+\frac{\e y}2)}\varepsilon+\frac{v\cdot v}{2\varepsilon}t \right)} \,
\overline{\psi}(x-\frac{\e y}2 + v t + x_0,t) 
e^{i\left(\frac{v\cdot (x-\frac{\e y}2)}\varepsilon+\frac{v\cdot v}{2\varepsilon}t \right)} dy=\\

=\int\limits_y e^{-2\pi i (k+\frac{v}{2\pi}) \cdot y} \psi(x+\frac{\e y}2 + v t + x_0,t) 
\overline{\psi}(x-\frac{\e y}2 + v t + x_0,t) 
 dy = W^\e[\psi(t)]\left(x+vt+x_0,k+\frac{v}{2\pi}\right).
\ea
\]

%
%
\qed

Lemma \ref{lm:GlInvII} holds for any Galilean invariant nonlinearity essentially with the same proof; in particular we have

\begin{lemma}[Galilean invariance] \label{lm:GlInv} Let $\psi$ satisfy
\beq \label{eq:GLIoreq}
i\varepsilon\partial_t \psi + \frac{\varepsilon^2}2 \Delta \psi - b|\psi|^{2\sigma} \psi=0, \qquad \psi(t=0)=\psi_0 \in L^2(\R^n).
\ee
For any $x_0,v\in \R^n$,
and denote
\beq \label{eq:GalTransf}
u(x,t) = \psi(x+vt+x_0,t) e^{-i\left(\frac{v\cdot x}\varepsilon+\frac{v\cdot v}{2\varepsilon}t \right)}.
\ee
Then $u$ satisfies
\beq \label{eq:aljh}
i\varepsilon\partial_t u + \frac{\varepsilon^2}2 \Delta u - b|u|^{2\sigma} u=0, \qquad u(t=0)=u_0=\psi_0(x+x_0)e^{-i\frac{v\cdot x}{\varepsilon}} \in L^2(\R^n).
\ee
Moreover,
\beq \label{eq:WTGalTransf}
W^\e[u(t)](x,k) = W^\e[\psi(t)]\left(x+vt+x_0,k+\frac{v}{2\pi}\right).
\ee
\end{lemma}

\begin{lemma}[Center of mass and conservation of momentum for the $1$-dimensional Schr\"odinger-Poisson] \label{lm:comsp}
Let $\psi$ satisfy
\beq \label{eq:GLIoreq22}
i\varepsilon\partial_t \psi + \frac{\varepsilon^2}2 \Delta \psi - \frac{b}2 \int\limits_x |x-y|  |\psi(y,t)|^{2} dy \, \psi=0, \qquad \psi(t=0)=\psi_0 \in \mathcal{S}(\R^n)
\ee

Then
\[
\frac{d}{dt} \mu_x(\psi(t)) = 2\pi \mu_k(\psi_0), \qquad  \frac{d}{dt} \mu_k(\psi(t))=0.
\]

\end{lemma}

\noindent 
{\bf Proof:} We compute
\[
\bac
\frac{d}{dt} \mu_x(\psi)=
\frac{d}{dt} \langle x \psi, \psi \rangle = \frac{i\e}{2} \left( \langle x\Delta \psi, \psi \rangle -  \langle x \psi, \Delta\psi \rangle \right)=\\
=
\frac{i\e}{2} \left( \langle  \psi, \nabla \psi \rangle -  \langle \nabla \psi, \psi \rangle \right)= i \e  \langle \nabla \psi, \psi \rangle = 2\pi \mu_k(\psi(t)).
\ea
\]
Moreover, denoting $V(x,t)=\frac{b}2 \int\limits_y |x-y|  |\psi(y,t)|^{2} dy$ the nonlinear potential we have
\[
\bac
\frac{d}{dt} \mu_k(\psi)=
\e \frac{d}{dt} \langle k \widehat \psi, \widehat \psi \rangle = 
\frac{\e}{2\pi i} \frac{d}{dt} \langle \nabla  \psi, \psi \rangle 
=\frac{1}\pi \real \langle \nabla \psi, V \psi \rangle = \\

= \frac{1}{2\pi} \int\limits_x V(x,t) (\overline{\psi} \nabla \psi+ \psi \nabla \overline{\psi}) dx = \frac{1}{2\pi} \int\limits_x V(x,t) \nabla |\psi(x,t)|^2 dx
\ea
\]
and now we complete the computation by observing 
\[
\bac
\int\limits_x V(x,t) \nabla |\psi(x,t)|^2 dx = \frac{b}2 \int\limits_x \int\limits_y |x-y|  |\psi(y,t)|^{2} dy \nabla |\psi(x,t)|^2 dx =\\

= -\frac{b}2 \int\limits_x \int\limits_y sign(x-y|)  |\psi(y,t)|^{2} dy  |\psi(x,t)|^2 dx = 0.
\ea
\]
\qed

\begin{lemma}[Center of mass and conservation of momentum for power nonlinearities] \label{lm:compn}
Let $\psi$ satisfy
\beq \label{eq:GLIoreqcom}
i\varepsilon\partial_t \psi + \frac{\varepsilon^2}2 \Delta \psi - b |\psi|^{2\sigma} \psi=0, \qquad \psi(t=0)=\psi_0 \in \mathcal{S}(\R^n)
\ee

Then
\[
\frac{d}{dt} \mu_x(\psi(t)) = 2\pi \mu_k(\psi_0), \qquad  \frac{d}{dt} \mu_k(\psi(t))=0.
\]

\end{lemma}

\noindent 
{\bf Proof:} We compute
\[
\bac
\frac{d}{dt} \mu_x(\psi)=
\frac{d}{dt} \langle x \psi, \psi \rangle = \frac{i\e}{2} \left( \langle x\Delta \psi, \psi \rangle -  \langle x \psi, \Delta\psi \rangle \right)=\\
=
\frac{i\e}{2} \left( \langle  \psi, \nabla \psi \rangle -  \langle \nabla \psi, \psi \rangle \right)= i \e  \langle \nabla \psi, \psi \rangle = 2\pi \mu_k(\psi(t)).
\ea
\]
Moreover, denoting $V(x,t)=b|\psi(x,t)|^{2\sigma}$ the nonlinear potential we have
\[
\bac
\frac{d}{dt} \mu_k(\psi)=
\e \frac{d}{dt} \langle k \widehat \psi, \widehat \psi \rangle = 
\frac{\e}{2\pi i} \frac{d}{dt} \langle \nabla  \psi, \psi \rangle 
=\frac{1}\pi \real \langle \nabla \psi, V \psi \rangle = \\

= \frac{1}{2\pi} \int\limits_x V(x,t) (\overline{\psi} \nabla \psi+ \psi \nabla \overline{\psi}) dx = \frac{1}{2\pi} \int\limits_x V(x,t) \nabla |\psi(x,t)|^2 dx
\ea
\]
and now we complete the computation by observing 
\[
\int\limits_x V(x,t) \nabla |\psi(x,t)|^2 dx = b \int\limits_x |\psi(x,t)|^{2\sigma} \nabla |\psi(x,t)|^2 dx = \frac{b}{\sigma+1} \int\limits_x \nabla |\psi(x,t)|^{2\sigma + 2} dx=0.
\]
\qed

\subsection{Inequalities}

\begin{lemma}\label{lm:bckr1} For any $a,b,q>0$
\[
\left( a+b \right)^q \leqslant C (a^q+b^q) \qquad \mbox{ for } \qquad C=\left\{ 
\begin{array}{l l}
2^{q-1}, & q \geqslant 1, \\
1, &  0 < q \leqslant 1
\ea
\right.
\]
\end{lemma}

\noindent {\bf Proof:} For $q\geqslant 1$, we use the convexity of $f(r)=r^q$, namely
\[
 f\left( \frac{a+b}{2}\right) \leqslant \frac{f(a)+f(b)}{2} \quad \Rightarrow \quad \left( \frac{a+b}{2}\right)^q \leqslant \frac{a^q+b^q}{2} \qquad \forall a,b >0
\]
For $q<1$,  $f(r)=r^q$ is concave and therefore sub-additive.
\qed

\begin{lemma}[Bootstrap argument] \label{lm:bootsrmine}
Let $f(t) \in C([0,\infty),[0,\infty))$, $0<A,B$, $0<\theta<1$ and
\[
f(t) \leqslant A + B f^\theta(t).
\]
Then $f(t)$ is bounded by the largest positive solution of
\beq \label{eq:lmbottstrehgf}
x- Bx^\theta -A =0.
\ee

In the case $\theta=\frac{1}2$,
\[
f(t) \leqslant A +\frac{B^2}2 + \frac{B \sqrt{B^2+4A}}2.
\]
\end{lemma}

\noindent {\bf Proof:}
Since $b\sqrt{t}$ grows more slowly than $t$ when $t\to \infty$, it is clear that $f(t)$ is bounded above.

Moreover the maximum value $f_{max}$ will satisfy \eqref{eq:lmbottstrehgf};
indeed if for some value $f$
\[
f<  A + B \sqrt{f}
\]
this means that a somewhat larger value $f$ would still be possible.

Thus we need to compute the largest solution of \eqref{eq:lmbottstrehgf}; if $\theta=\frac{1}2$ this is achieved by solving the quadratic equation 
\[
\left(\sqrt{f_{max}}\right)^2 - B \sqrt{f_{max}} - A =0.
\]
\qed

\begin{theorem}[Gagliardo-Nirenberg $L^2$-gradient inequality] \label{thrm:GLineq} For every $f$ such that if
\[
f \in L^q, \qquad \nabla f \in L^2
\]
 then
\[
\|f\|_{L^p(\R^n)} \leqslant C^{GN}_{q,p,n} \|\nabla f \|_{L^2(\R^n)}^\theta  \|f\|_{L^q(\R^n)}^{1-\theta}
\]
for
\[
\bal
& 1 < q < p < \frac{2n}{(n-2)_+}  
\eal
\]
and
\[
\theta = 
\frac{2n(p-q)}{p[2n-q(n-2)]}.
\]
Moreover, the sharp constant $C^{GN}_{q,p,n}$ is known.
\end{theorem}

\noindent {\bf Proof:} See \cite{agueh}.
\qed

\begin{corollary}\label{cor:GN22} Let 
\[
f \in H^1(\R^n), \qquad 
  \|f\|_{L^2(\R^n)}=1, \qquad
\sigma \in \left(0,\frac{2}{(n-2)_+}\right).
\]
Then 
\[
\|f\|^{2\sigma+2}_{L^{2\sigma+2}(\R^n)} \leqslant C^{GN}_*\, \|\nabla f \|_{L^2(\R^n)}^{{n\sigma}}.
\]
\end{corollary}

\noindent {\bf Proof:} Set $q=2$ and  $p=2\sigma+2$ in Theorem \ref{thrm:GLineq}. 
The constant $C^{GN}_*(n,\sigma):=\big( C^{GN}_{2,2\sigma+2,n} \big)^{2\sigma+2}$ is sharp and known \cite{agueh}.
\qed

\subsection{Computations for squeezed states}

\begin{definition} \label{def:sqstates} Let 
\[
a\in \mathcal{S}(\R^n), \qquad \|a\|_{L^2}=1, \qquad \mu_x(a) = \mu_k(a)=0,
\]
$\beta \in (0,1)$. The function
\[
\psi^\e_0(x) = \e^{-\frac{n\beta}2} a(\frac{x-x_0}{\e^\beta}) e^{\frac{2\pi i k_0 \cdot (x-x_0)}\e}
\] 
will be called a squeezed state with envelope $a$ and rate of concentration $\beta$.
\end{definition}

\begin{lemma}Let 
\[
\psi^\e_0(x) = \e^{-\frac{n\beta}2} a(\frac{x-x_0}{\e^\beta}) e^{\frac{2\pi i k_0 \cdot (x-x_0)}\e}
\] 
be  a squeezed state with envelope $a$ and rate of concentration $\beta$. Then
\[
\|\psi^\e_0\|_{L^2}=1, \qquad 
\mu_x(\psi_0^\e) = x_0, \qquad \mu_k(\psi_0^\e) = k_0, \qquad \sigma_x(\psi_0^\e)= O(\e^{\beta})
\qquad \sigma_k(\psi_0^\e)= O(\e^{1-\beta}).
\]
\end{lemma}

\noindent {\bf Proof:}
One readily computes
\[
\sigma_x(\psi_0^\e)= \| x \e^{-\frac{n\beta}2} a(\frac{x}{\e^\beta}) \|_{L^2}= O(\e^{{\beta}}),
\]
and
\[
\sigma_k(\psi_0^\e)= \| \e^{1-\frac{n\beta}2} \nabla  a(\frac{x}{\e^\beta}) \|_{L^2}= O(\e^{1-{\beta}}).
\]
\qed

\section{Proof of the main results} \label{sec:mainproofs}

\subsection{Proof of Theorem \ref{thrm:Main1}} \label{subsec:mainproofs1}

By virtue of Lemma \ref{lm:GlInvII}, 
the solution of the problem
\beq \label{eq:peioygg}
i\varepsilon\partial_t u^\varepsilon + \frac{\varepsilon^2}2 \Delta u^\varepsilon - \frac{b}2 \int\limits_y |x-y| |u^\varepsilon(y,t)|^{2}dy \, u^\varepsilon =0 \qquad u^\e(t=0)=u^\e_0=\psi^\e_0(x+\mu_x(\psi^\e_0))e^{-2\pi i\frac{\mu_k(\psi^\e_0)\cdot x}{\varepsilon}} 
\ee
is related to $\psi^\e$ through 
\beq 
u^\e(x,t) = \psi^\e(x+vt+x_0,t) e^{-i\left(\frac{v\cdot x}\varepsilon+\frac{v\cdot v}{2\varepsilon}t \right)}, \qquad v=2\pi \mu_k(\psi^\e_0), \quad  x_0=\mu_x(\psi^\e_0).
\ee
By virtue of Lemma \ref{lm:comsp} and by the construction of $u_0^\e$,
\beq \label{eq:onmoe66y}
\bac
\mu_x(u^\e(t))=\mu_k(u^\e(t)) = 0, \\ \sigma_x(\psi^\e(t)) = \sigma_x(u^\e(t)) = \| x u^\e(t) \|_{L^2_x}, \qquad 
\sigma_k(\psi^\e(t)) = \sigma_k(u^\e(t)) = \frac{1}{2\pi} \| \e \nabla u^\e(t)   \|_{L^2}.
\ea
\ee

For now we will work with equation \eqref{eq:peioygg}, and ultimately transfer our results to $W^\e[\psi^\e(t)]$.

By virtue of the conservation of energy \eqref{eq:spconoe}, we have
\[
\bac
\frac{\e^2}{2} \|\nabla u^\e(t)\|^2_{L^2} \leqslant 
\frac{\varepsilon^2}2 \|\nabla u^\varepsilon(t)\|_{L^2}^2 + \frac{b}{4} \int\limits_{x,y} |x-y| |u^\e(x,t)|^2 |u^\e(y,t)|^2 dxdy = \\

=\frac{\varepsilon^2}2 \|\nabla u^\varepsilon_0\|_{L^2}^2 + \frac{b}{4} \int\limits_{x,y} |x-y| |u^\e_0(x)|^2 |u^\e_0(y)|^2 dxdy \leqslant 

\frac{\varepsilon^2}2 \|\nabla u^\varepsilon_0\|_{L^2}^2 + \frac{b}2 \int\limits_{x} |x| |u^\e_0(x)|^2  dx \leqslant\\

\leqslant\frac{\varepsilon^2}2 \|\nabla u^\varepsilon_0\|_{L^2}^2 + \frac{b}2 \| x u_0^\e\|_{L^2}.
\ea
\]
Thus by virtue of Lemma \ref{lm:bckr1} we have
\beq \label{eq:hgdryu}
\| \e \nabla u^\e(t)\|_{L^2} \leqslant \| \e \nabla u^\e_0\|_{L^2} + \sqrt{b \|x u^\e_0\|_{L^2} }.
\ee
Moreover, by virtue of equation \eqref{eq:momscrpois},
\beq \label{eq:hgdryu12}
\| x u^\e(t)\|_{L^2} \leqslant \|x u^\e_0\|_{L^2} + t \left( \| \e \nabla u^\e_0\|_{L^2} + \sqrt{b \|x u^\e_0\|_{L^2}} \right).
\ee
Recalling equation \eqref{eq:onmoe66y}, we can recast equations \eqref{eq:hgdryu}, \eqref{eq:hgdryu12} as
\[
\bac 
\sigma_k(\psi^\e(t))=
\sigma_k(u^\e(t)) \leqslant \sigma_k(u^\e_0) + \sqrt{\frac{b}{2\pi} \sigma_x(u^\e_0)}, \qquad 
\sigma_x(\psi^\e(t))=
\sigma_x(u^\e(t)) \leqslant \sigma_x(u^\e_0) + t \left( \sigma_k(u^\e_0) + \sqrt{\frac{b}{2\pi} \sigma_x(u^\e_0)} \right),
\ea
\]
and finally 
\[
\sigma_k(\psi^\e(t)) + \sigma_x(\psi^\e(t)) \leqslant \sigma_k(u^\e_0) (1+t) + 
\sigma_x(u^\e_0) + \sqrt{\frac{b}{2\pi}}(1+t) \sqrt{\sigma_x(u_0^\e)}.
\]
The proof is complete by recalling that
\[
 \mu_x(\psi^\e(t)) =  \mu_x(\psi^\e_0) + 2\pi t \mu_k(\psi^\e_0), \qquad   \mu_k(\psi^\e(t))=\mu_k(\psi^\e_0),
\]
by virtue of Lemma \ref{lm:comsp}, and then applying Corollary \ref{cor:wavpkts}.
\qed 

\subsection{Proof of Theorem \ref{thrm:Main2}} \label{subsec:mainproofs2}

In exact analogy to what we did before,
the solution of the problem
\beq \label{eq:peioygg20}
i\varepsilon\partial_t u^\varepsilon + \frac{\varepsilon^2}2 \Delta u^\varepsilon - b |u^\e|^{2\sigma} \, u^\varepsilon =0 \qquad u^\e(t=0)=u^\e_0=\psi^\e_0(x+\mu_x(\psi^\e_0))e^{-2\pi i\frac{\mu_k(\psi^\e_0)\cdot x}{\varepsilon}} 
\ee
is related to $\psi^\e$ through 
\beq 
u^\e(x,t) = \psi^\e(x+vt+x_0,t) e^{-i\left(\frac{v\cdot x}\varepsilon+\frac{v\cdot v}{2\varepsilon}t \right)}, \qquad v=2\pi \mu_k(\psi^\e_0), \quad  x_0=\mu_x(\psi^\e_0).
\ee
Again, by virtue of Lemma \ref{lm:comsp} and by the construction of $u_0^\e$,
\beq \label{eq:on0oe663y}
\bac
\mu_x(u^\e(t))=\mu_k(u^\e(t)) = 0, \\ \sigma_x(\psi^\e(t)) = \sigma_x(u^\e(t)) = \| x u^\e(t) \|_{L^2_x}, \qquad 
\sigma_k(\psi^\e(t)) = \sigma_k(u^\e(t)) = \frac{1}{2\pi} \| \e \nabla u^\e(t)   \|_{L^2}.
\ea
\ee

By virtue of the conservation of energy, equation \eqref{eq:defocenpn},
\[
\bac
\frac{\e^2}2\|\nabla u^\e(t)\|_{L^2}^2 \leqslant \frac{\e^2}2\|\nabla u^\e_0\|_{L^2}^2 + \frac{b}{\sigma +1} \|u^\e_0\|_{L^{2\sigma+2}}^{2\sigma+2} \leqslant 
\frac{1}2\|\e\nabla u^\e_0\|_{L^2}^2 + \frac{b}{\sigma +1} C^{GN}_* \|\nabla u^\e_0\|_{L^2}^{n\sigma}, 
%
\ea
\]
where in the last step we used the Gagliardo-Nirenberg inequality, Corollary \ref{cor:GN22}. Using Lemma \ref{lm:bckr1}, this becomes
\[
\|\e\nabla u^\e(t)\|_{L^2}  \leqslant \|\e\nabla u^\e_0\|_{L^2} \,  + \, \sqrt{ \e^{-n\sigma} \, b \,\, \frac{C^{GN}_*}{2\sigma +2} \,\,} \, \|\e\nabla u^\e_0\|_{L^2}^{\frac{n\sigma}2}
\]
Moreover, equation \eqref{eq:momscrpois22} of Theorem \ref{thrm:momscrpois22} implies that
\[
\|xu^\e(t)\|_{L^2} \leqslant \|xu^\e_0\|_{L^2} + t \left( \|\e\nabla u^\e_0\|_{L^2}  + \sqrt{ \e^{-n\sigma} \, b \,\, \frac{C^{GN}_*}{2\sigma +2} \,\,} \, \|\e\nabla u^\e_0\|_{L^2}^{\frac{n\sigma}2}
\right).
\]
Collecting the last two equations, and recalling equation \eqref{eq:on0oe663y}, we have
\[
\sigma_x(u^\e(t)) + \sigma_k(u^\e(t)) \leqslant 
\sigma_x(u^\e_0) + ( 1 + t ) \sigma_k(u^\e_0)   + (1+t)\left( \sigma_k^{\frac{n\sigma}2}(u_0^\e) \, \sqrt{  \frac{b}{\e^{n\sigma}} \,\, \frac{C^{GN}_* (2\pi)^{n\sigma-2}}{2\sigma +2} \,\,} 
  \right)
\]
The proof is complete by recalling that
\[
 \mu_x(\psi^\e(t)) =  \mu_x(\psi^\e_0) + 2\pi t \mu_k(\psi^\e_0), \qquad   \mu_k(\psi^\e(t))=\mu_k(\psi^\e_0),
\]
by virtue of Lemma \ref{lm:compn}, and then applying Corollary \ref{cor:wavpkts}.
\qed

\subsection{Proof of Theorem \ref{thrm:Main3}} \label{subsec:mainproofs3}

In exact analogy to what we did before,
the solution of the problem
\beq \label{eq:peijghoygg20}
i\varepsilon\partial_t u^\varepsilon + \frac{\varepsilon^2}2 \Delta u^\varepsilon - b |u^\e|^{2\sigma} \, u^\varepsilon =0 \qquad u^\e(t=0)=u^\e_0=\psi^\e_0(x+\mu_x(\psi^\e_0))e^{-2\pi i\frac{\mu_k(\psi^\e_0)\cdot x}{\varepsilon}} 
\ee
is related to $\psi^\e$ through 
\beq 
u^\e(x,t) = \psi^\e(x+vt+x_0,t) e^{-i\left(\frac{v\cdot x}\varepsilon+\frac{v\cdot v}{2\varepsilon}t \right)}, \qquad v=2\pi \mu_k(\psi^\e_0), \quad  x_0=\mu_x(\psi^\e_0).
\ee
Again, by virtue of Lemma \ref{lm:comsp} and by the construction of $u_0^\e$,
\beq \label{eq:on0oe6hgjc63y}
\bac
\mu_x(u^\e(t))=\mu_k(u^\e(t)) = 0, \\ \sigma_x(\psi^\e(t)) = \sigma_x(u^\e(t)) = \| x u^\e(t) \|_{L^2_x}, \qquad 
\sigma_k(\psi^\e(t)) = \sigma_k(u^\e(t)) = \frac{1}{2\pi} \| \e \nabla u^\e(t)   \|_{L^2}.
\ea
\ee

By virtue of the conservation of energy, equation \eqref{eq:defocenpn},
\[
\bac
\frac{\e^2}2\|\nabla u^\e(t)\|_{L^2}^2 = \frac{\e^2}2\|\nabla u^\e_0\|_{L^2}^2 + \frac{b}{\sigma +1} \|u^\e_0\|_{L^{2\sigma+2}}^{2\sigma+2} + \frac{|b|}{\sigma +1} \|u^\e(t)\|_{L^{2\sigma+2}}^{2\sigma+2} \leqslant \\

\leqslant \frac{\e^2}2\|\nabla u^\e_0\|_{L^2}^2 + \frac{|b|}{\sigma +1} \|u^\e(t)\|_{L^{2\sigma+2}}^{2\sigma+2} \leqslant

\frac{1}2\|\e\nabla u^\e_0\|_{L^2}^2 + \frac{|b|}{\sigma +1} C^{GN}_* \|\nabla u^\e(t)\|_{L^2}^{n\sigma}, 
%
\ea
\]
where in the last step we used the Gagliardo-Nirenberg inequality, Corollary \ref{cor:GN22}. Using Lemma \ref{lm:bckr1}, this becomes
\[
\|\e\nabla u^\e(t)\|_{L^2}  \leqslant \|\e\nabla u^\e_0\|_{L^2} +  \sqrt{ \e^{-n\sigma} \, |b(\e)| \,\, \frac{C^{GN}_*}{2\sigma +2} \,\,}  \|\e\nabla u^\e(t)\|_{L^2}^{\frac{n\sigma}2}.
\]
Since $\frac{n\sigma}2=\frac{1}2$, Lemma \ref{lm:bootsrmine} applies to $f(t)=\|\e\nabla u^\e(t)\|_{L^2}$, yielding
\beq \label{eq:fockinetcibd}
\|\e\nabla u^\e(t)\|_{L^2}  \leqslant \|\e\nabla u^\e_0\|_{L^2} +\frac{|b(\e)| C^{GN}_*}{\e(4+\frac{4}n)} + \frac{1}{2} \sqrt{   \frac{|b(\e)| C^{GN}_*}{ \e (2 +\frac{2}n)} \,\,}
 \sqrt{    \frac{|b(\e)| C^{GN}_*}{ \e (2 +\frac{2}n)}  + 4\|\e\nabla u^\e_0\|_{L^2}
 }
\ee
For brevity we will denote
\beq \label{eq:kinbdsymbol}
\mathcal{K} := \frac{|b(\e)| C^{GN}_*}{\e(4+\frac{4}n)} + \frac{1}{2} \sqrt{   \frac{|b(\e)| C^{GN}_*}{ \e (2 +\frac{2}n)} \,\,}
 \sqrt{    \frac{|b(\e)| C^{GN}_*}{ \e (2 +\frac{2}n)}  + 4\|\e\nabla u^\e_0\|_{L^2}
 }
\ee

Moreover, equation \eqref{eq:momscrpois22} of Theorem \ref{thrm:momscrpois22} implies that
\[
\|xu^\e(t)\|_{L^2} \leqslant \|xu^\e_0\|_{L^2} + t \Big(\|\e\nabla u^\e_0\|_{L^2}  + \mathcal{K} \Big).
\]
Collecting the last two equations, and recalling equation \eqref{eq:on0oe6hgjc63y}, we have
\[
\sigma_x(u^\e(t)) + \sigma_k(u^\e(t)) \leqslant 
\sigma_x(u^\e_0) + \sigma_k(u^\e_0) ( 1 + t) +  \mathcal{K} \frac{1+t}{2\pi}.
\]

The proof is complete by recalling that
\[
 \mu_x(\psi^\e(t)) =  \mu_x(\psi^\e_0) + 2\pi t \mu_k(\psi^\e_0), \qquad   \mu_k(\psi^\e(t))=\mu_k(\psi^\e_0),
\]
by virtue of Lemma \ref{lm:compn}, and then applying Corollary \ref{cor:wavpkts}.
\qed

\subsection{The concentration estimates}

\begin{lemma}[Concentration of Wigner transforms to $\delta(x-0,k-0)$ for Schwarz functions]\label{lm:loc11} Let $u \in \mathcal{S}(\R^n)$. Then
\[
\|  {W}^\varepsilon[u]- \|u\|_{L^2}^2 \cdot \delta(x-0,k-0)\|_{A^{-1}}  \leqslant 
\|u\|_{L^2} \left( 2 \pi  \|x u \|_{L^2} + \varepsilon \|\nabla u \|_{L^2}  \right).
\]
\end{lemma}

\noindent {\bf Proof: } For brevity  we will denote $W^\varepsilon(x,k)=W^\varepsilon[u](x,k)$, and $X,K$ the Fourier dual variables to $x,k$. Naturally, the idea of the proof will be to work on the Fourier dual of the variables in which the Lemma is stated, namely we will use the fact that
\[
\left| \langle {W}^\varepsilon- \|u\|_{L^2}^2 \cdot \delta(x-0,k-0), {\phi} \rangle \right| =\left| \langle \widehat{W}^\varepsilon(X,K)-\|u\|_{L^2}^2, \widehat{\phi} \rangle \right|.
\]
In what follows we will use the elementary computation
\beq \label{eq:wignFour}
\widehat{W}^\varepsilon(X,K)=\mathcal{F}_{(x,k) \to (X,K)} [W^\e(x,k)] = \int\limits_{x} e^{-2\pi i x \cdot X} u(x-\frac{\varepsilon K}2) \overline{u(x+\frac{\varepsilon K}2)} dx.
\ee

Now observe that, for any $j\in \{ 1,\dots,n\}$,
\beq \label{wq:tricktK1}
\bac
\partial_{K_j} \widehat{W}^\e(X,K)=
\partial_{K_j} \int\limits_{x} e^{-2\pi i xX} u(x-\frac{\varepsilon K}2) \overline{u(x+\frac{\varepsilon K}2)} dx =\\
 = 
 \frac{\varepsilon}2 \int\limits_{x} e^{-2\pi i xX} \left[{ 
u(x-\frac{\varepsilon K}2) \partial_{x_j} \overline{u(x+\frac{\varepsilon K}2)} - \overline{u(x+\frac{\varepsilon K}2)} \partial_{x_j} u(x-\frac{\varepsilon K}2)
}\right] dx
  \quad \Rightarrow\\
  
  \Rightarrow \quad  |\partial_{K_j} \widehat{W}^\varepsilon(X,K)| \leqslant \varepsilon \|\nabla u \|_{L^2} \|u\|_{L^2},
\ea
\ee
where we used the fact that
\[
\left| \int\limits_{x} e^{-2\pi i x \cdot X} u(x-\frac{\varepsilon K}2) \overline{v(x+\frac{\varepsilon K}2)} dx \right| \leqslant \|u\|_{L^2}\|v\|_{L^2}
\]
by virtue of the Cauchy-Schwarz inequality.

On the other hand,
\beq
\bac \label{eq:xderwhat}
\frac{i}{\pi} \partial_{X_j} \widehat{W}^\e(X,K) = \mathcal{F}_{(x,k)\to(X,K)}[2 x_j W^\e(x,k)]=\\
=2\int\limits_{x,k} e^{-2\pi i [kK+xX] }  x_j  W^\varepsilon(x,k) dxdk 
=
2 \int\limits_{x,k,y} e^{-2\pi i k[K+y] }dk   e^{-2\pi i xX }  x_j  u(x+\frac{\varepsilon y}2) \overline{u(x-\frac{\varepsilon y}2)} dxdy =\\

= 2 \int\limits_{x}    e^{-2\pi i xX }  x_j  u(x-\frac{\varepsilon K}2) \overline{u(x+\frac{\varepsilon K}2)} dx=\\
=  \int\limits_{x}    e^{-2\pi i xX }   \left[{ 
(x-\frac{\varepsilon K}2)  u(x-\frac{\varepsilon K}2) \overline{u(x+\frac{\varepsilon K}2)}  +   u(x-\frac{\varepsilon K}2) (x+\frac{\varepsilon K}2) \overline{u(x+\frac{\varepsilon K}2)}
 }\right] 
 dx\quad \Rightarrow \\
 
 \Rightarrow \quad |\partial_{X_j} \widehat{W}^\varepsilon(X,K)| \leqslant 2 \pi \| u\|_{L^2} \|x_j u \|_{L^2}.
\ea
\ee

Combining equations \eqref{wq:tricktK1} and \eqref{eq:xderwhat} it follows that
\beq
\| \nabla_{X,K} \widehat{W}(X,K) \|_{L^\infty_{X,K}} \leqslant \|u\|_{L^2} \left( 2\pi\|x u\|_{L^2} + \e \|\nabla_x u \|_{L^2} \right).
\ee

Finally, observe that 
\beq
\widehat{W}^\e(0,0)=\|u\|_{L^2}^2,
\ee
e.g. by evaluating equation \eqref{eq:wignFour} at $(X,K)=(0,0)$. Now we Taylor expand $\widehat{W}(X,K)$ around $(0,0)$ to obtain
\beq \label{eq:wmtexp}
\bac
\left|\widehat{W}^\e(X,K) - \|u\|_{L^2}^2\right| \leqslant |(X,K)| \cdot \|\nabla_{X,K} \widehat{W} \|_{L^\infty} 

 \leqslant \sqrt{|X|^2+|K|^2} \,\, \|u\|_{L^2} \left( 2 \pi  \|x u \|_{L^2} + \varepsilon \|\nabla u \|_{L^2}  \right).
\ea
\ee

The proof is completed by integrating against any $A^1$ test function $\phi$,
\beq
\bac
\left| \langle {W}^\varepsilon- \|u\|_{L^2}^2 \cdot \delta(0,0), {\phi} \rangle \right| =\left| \langle \widehat{W}^\varepsilon(X,K)-\|u\|_{L^2}^2, \widehat{\phi} \rangle \right| \leqslant \\

\leqslant
\|u\|_{L^2} \left( 2 \pi  \|x u \|_{L^2} + \varepsilon \|\nabla u \|_{L^2}  \right) \int\limits_{X,K} \sqrt{|X|^2+|K|^2}  |\widehat{\phi}(X,K)|dXdK\leqslant\\

\leqslant\|u\|_{L^2} \left( 2 \pi  \|x u \|_{L^2} + \varepsilon \|\nabla u \|_{L^2}  \right) \|\phi\|_{A^1}.
\ea
\ee

\qed

\begin{corollary}[Concentration of Wigner transforms to $\delta(\mu_x(\psi),\mu_k(\psi))$ for Sobolev functions] \label{cor:wavpkts} Let $\psi \in H^1 \cap \widehat{H}^1$, $\|\psi\|_{L^2}=1$. Then
\beq \label{eq:cor52a}
\|  {W}^\varepsilon[\psi]-  \delta(x-0,k-0)\|_{A^{-1}}  \leqslant 
 2 \pi  \| x \psi \|_{L^2} + \varepsilon \|\nabla \psi \|_{L^2},
\ee
and more generally
\beq \label{eq:cor52b}
\left\|  {W}^\varepsilon[\psi]-  \delta\big(x-\mu_x(\psi),k-\mu_k(\psi)\big)\right\|_{A^{-1}}  \leqslant 
2\pi \Big( \sigma_x(\psi) + \sigma_k(\psi) \Big).
\ee
\end{corollary}

\noindent {\bf Proof:} The proof of the Corollary consists of two parts:  first, we check that the arguments in the proof of Lemma \ref{lm:loc11} still work for $H^1 \cap \widehat{H}^1$ wavefunctions. Then we apply a Galilean transform to obtain  concentration on any point of phase-space.

Since $\psi\in H^1(\R^n)\cap \widehat{H}^1(\R^n)$ a basic computation shows that $W^\e[\psi] \in H^1(\R^{2n})\cap \widehat{H}^1(\R^{2n}) \cap L^\infty(\R^{2n})$ \cite{AP3,AMP}. Moreover, equations \eqref{wq:tricktK1} and \eqref{eq:xderwhat} mean that $\widehat{W}^\e[\psi] \in W^{1,\infty}(\R^{2n})$.
Therefore the Taylor expansion of equation \eqref{eq:wmtexp} makes sense as a Taylor expansion in $W^{1,\infty}(\R^{2n})$ \cite{LL}, and equation \eqref{eq:cor52a} follows.

In order to prove equation \eqref{eq:cor52b}, let us call $u$ the ``centered version of $\psi$,'' 
\[
u(x) = \mathcal{M}_{\frac{\mu_k(\psi)}{\e} }\mathcal{T}_{\mu_x(\psi)} \psi=\psi(x+\mu_x(\psi))e^{-2\pi i\frac{\mu_k(\psi)\cdot x}{\varepsilon}};
\]
by construction $\mu_x(u)=\mu_k(u)=0$.
Now observing that 
\[
\sigma_x(\psi) = \sigma_x(u) = \| x u \|_{L^2}, \qquad \sigma_k(\psi) = \sigma_k(u) = \frac{\varepsilon}{2\pi} \|\nabla u \|_{L^2},
\]
equation \eqref{eq:cor52a} implies that
\beq \label{eq:wtmacli12}
\|  {W}^\varepsilon[u]-  \delta(x-0,k-0)\|_{A^{-1}}  \leqslant 
2\pi \Big( \sigma_x(\psi) + \sigma_\mu(\psi) \Big).
\ee
Moreover, 
\[
\bac
W^\e[u(t)] = W^\e[\psi(x+x_0) e^{-i\frac{2\pi \mu_k(\psi)\cdot x}\varepsilon }] = \\

=\int\limits_y e^{-2\pi i k \cdot y} \psi(x+\frac{\e y}2  + \mu_x(\psi)) 
e^{-i\frac{2\pi \mu_k(\psi)\cdot (x+\frac{\e y}2)}\varepsilon} \,
\overline{\psi}(x-\frac{\e y}2  + \mu_x(\psi)) 
e^{i\frac{2\pi \mu_k(\psi)\cdot (x-\frac{\e y}2)}\varepsilon} dy=\\

=\int\limits_y e^{-2\pi i (k+\mu_k(\psi)) \cdot y} \psi(x+\frac{\e y}2  + \mu_x(\psi)) 
\overline{\psi}(x-\frac{\e y}2 +  \mu_x(\psi)) 
 dy = W^\e[\psi]\big(x+\mu_x(\psi),k+\mu_k(\psi)\big)
\ea
\]
and thus \eqref{eq:wtmacli12} means
\[
\bac
\|  W^\e[\psi]\left(x+\mu_x(\psi),k+\mu_k(\psi)\right)-  \delta(x-0,k-0)\|_{A^{-1}}  \leqslant 
2\pi \Big( \sigma_x(\psi) + \sigma_k(\psi)  \Big) \quad \Leftrightarrow \\

\Leftrightarrow \quad 
\|  W^\e[\psi](x,k)-  \delta\left(x-\mu_x(\psi),k-\mu_k(\psi)\right)\|_{A^{-1}}  \leqslant 
2\pi \Big( \sigma_x(\psi) + \sigma_k(\psi)  \Big).
\ea
\]
\qed

{\small

\end{document}